\documentclass[11pt,twoside]{article} 
\usepackage{amssymb,amsmath,amsfonts,amsthm}
\usepackage{amsmath}
\usepackage{graphicx,graphics}
\usepackage{latexsym}
\usepackage[applemac]{inputenc}
\usepackage{ae,aecompl}
\usepackage{mathrsfs}
\usepackage[english]{babel}
 \usepackage[colorlinks=true]{hyperref}

\renewcommand{\leq}{\lesqlant}

\newcommand{\Z}{\mathbb Z}

\newcommand{\ligne}{\vspace{3mm}}

\newsavebox{\fmbox}

\renewcommand{\geq}{\geqslant}
\renewcommand{\leq}{\leqslant}

\newcommand{\op}[1]{\operatorname{#1 }}

\newtheorem{theorem}{Theorem}

\newtheorem{lemma}[theorem]{Lemma}

\newtheorem{corollary}[theorem]{Corollary}
\newtheorem{proposition}[theorem]{Proposition}

\newtheorem{quest}[theorem]{Question}
\newtheorem{rek}{Remark}
\theoremstyle{definition}
\numberwithin{equation}{section}

\title{Recurrence of the $ \mathbb{Z}^d$-valued infinite snake \emph{via} unimodularity}

\author{Itai Benjamini and Nicolas Curien}
\date{}
\begin{document}
\maketitle

\begin{abstract}
We use the concept of unimodular random graph to show that the branching simple random walk on $\Z^{d}$ indexed by a critical geometric Galton-Watson tree conditioned to survive is recurrent if and only if $d \leq 4$.
\end{abstract}

\section*{Introduction}
Consider a simple random walk on $\mathbb{Z}^d$ indexed by some tree $T$. If $T$ is a super-critical Galton-Watson tree, the process we obtain is a branching random walk on $\mathbb{Z}^d$. It has been shown by Biggins \cite{Big77} that this branching random walk is almost surely recurrent \textit{i.e.}, visits the origin of $\mathbb{Z}^d$ infinitely often. When $T$ is a critical Galton-Watson tree, the walk is closely related to the theory of superBrownian motion and the associated  random snake of Le Gall, see \cite{Ald93b,LG99}. In this note we study the simple random walk on $\mathbb{Z}^d$ indexed by the critical geometric Galton-Waton tree $T_{\infty}$ \emph{conditioned to survive} \cite{Kes86}. Specifically we prove:
\begin{theorem} \label{main}
The simple random walk on $\Z^d$ indexed by $T_{\infty}$ is recurrent if and only if  $d \leq 4$.
\end{theorem} Recurrence of arbitrary Markov chains indexed by arbitrary trees was studied in \cite{BP94}, but the theorem above was not covered. Notice that the critical role of dimension $4$ is reminiscent of the theory of superBrownian motion (see \cite{LG99}) where the continuous analogue of Theorem \ref{main} is known \cite{SW94}.  

The proof of Theorem \ref{main} is based on the use of the ``Mass Transport Principle'' (see \cite[Section 3.2]{BS01}, and \cite{AL07}) and the related concept of unimodular random graphs combined with simple geometric estimates regarding the tree $T_{\infty}$. In particular we do not use any calculations related to the theory of superBrownian motion. More important than the application to Theorem \ref{main}, we believe that our technique could be applied in a much wider setup, see Section \ref{extensions}.

\bigskip The note is organized as follows. In the first section we recall the definition of the random  infinite tree $T_\infty$ and gather some simple geometric estimates. We also establish, using stationarity of the tree after re-rooting along a random walk, that the random graph $T_\infty$ satisfies the Mass Transport Principle. Section $2$ contains the proof of Theorem \ref{main}. We end the note with a few extensions, comments and open questions.
\bigskip

\noindent \textbf{Acknowledgments:} We are indebted to Ofer Zeitouni for a useful discussion. Thanks also go to Thomas Duquesne, Jean-François Le Gall, Yuval Peres and to an anonymous referee for valuable comments.

  \section{Definition and properties of $T_\infty$}

\subsection{Uniform plane trees}
A rooted  tree $\tau$ is a tree in the graph theoretic sense with a distinguished vertex $\rho$ called the root vertex. The tree $\tau$ can thus be seen as a family tree with ancestor $\rho$. A rooted ordered tree (or plane tree) is a rooted tree for which we have specified an ordering for the children of each vertex. See \cite{LG06} for a detailed definition. If $u$ is a vertex of $\tau$, we denote by $ \mathrm{deg}(u)$ the \emph{degree} of $u$ in $\tau$, that is its number of incident edges.

For $n \geq 1$, we denote  the set of all  plane trees with $n$ edges by $\mathcal{T}_{n}$.  For convenience,  we  associate with each tree $\tau \in \mathcal{T}_{n}$ a distinguished oriented edge $\vec{e}$ going from the first child of $\rho$ towards $\rho$. In the following $T_{n}$ is a random variable uniformly distributed over $\mathcal{T}_{n}$ and conditionally on $T_{n}$,  $X_{1}$ is a one-step simple random walk on $T_{n}$ starting at $\rho$: Equivalently $X_1$ is a uniform neighbor of the root $\rho$. We also keep track of the transition by indicating the  oriented edge $(\rho,X_1)$. We denote by $T_{n}^{(1)}$ the plane tree obtained from the  tree $T_{n}$ by keeping the planar ordering and by changing the distinguished oriented edge from $ \vec{e}$ to $(\rho,X_{1})$.

\begin{proposition} \label{prop:rev}The random plane tree $T_{n}^{(1)}$ is uniformly distributed over $\mathcal{T}_{n}$.

\end{proposition}
\proof For every given tree $\tau \in \mathcal{T}_{n}$ with oriented edge $ \vec{e}=(e_{-},e_{+})$, it is easy to see that exactly $\mathrm{deg}(e_{-})$ plane trees can give rise to $\tau$ after changing the oriented edge by a one step random walk: They consist of all the trees obtained from $\tau$ after exchanging  the oriented edge with an edge targeting $e_{-}$. Each of these trees has a probability $1/\mathrm{deg}(e_{-})$ to be transformed into $\tau$ after a one-step simple random walk. Thus $T_{n}^{(1)}$ is a uniform plane tree with $n$ edges. 
\endproof

If we denote $E_{-}$ and $E_{+}=\rho$ the origin and target vertices of the distinguished edge of $T_n$, we deduce that  $ (T_n,E_-,E_+)$  and $(T_n^{(1)},\rho,X_1)$ have the same distribution as random graphs with two distinguished neighboring vertices. Since the law of $T_n$ is unchanged under reversion of the distinguished edge, we get  \begin{eqnarray} \label{reversible} (T_n,\rho,X_1) & \overset{(d)}{=}& (T_n,X_1,\rho).  \end{eqnarray} In other words, the random rooted graph obtained from $T_n$ after forgetting the planar structure is \emph{reversible} in the sense of \cite[Definition 1]{BC10b}.

\subsection{The uniform infinite plane tree}

If $\tau$ is a plane tree and $k \in \{0,1,2, \ldots \}$, we denote by $[\tau]_{k}$ the plane tree obtained from $\tau$ by keeping the first $k$ generations of $\tau$ only. Let $T$ be a Galton-Watson (plane) tree with geometric offspring distribution of parameter $1/2$. It is classical that the distribution of $T$ conditionally on $T$ having $n$ edges is uniform over $\mathcal{T}_{n}$, see \cite{LG06}. Using this observation, the result of \cite{Dur03} (which has been folklore for a long time, see \cite{Ald91c,Kes86} and \cite{Gri80} for the case of unordered trees) can be interpreted as follows:
\begin{lemma} \label{convergence} Let $T_{n}$ be uniformly distributed over $\mathcal{T}_{n}$. Then there exists a random infinite plane tree $T_{\infty}$ such that for every $k\geq 0$ we have 
\begin{eqnarray} 
[T_{n}]_{k} & \xrightarrow[n\to \infty]{(d)} & [T_{\infty}]_{k}.
\end{eqnarray}
The random infinite plane tree $T_{\infty}$ is called \emph{the uniform infinite plane tree} or \emph{the geometric Galton-Watson tree conditioned to survive}.
\end{lemma}

The random infinite plane tree $T_{\infty}$ can be described as follows: Start with a semi-infinite line of vertices called the spine of the tree and graft to the left and to the right of each vertex of the spine an independent  critical geometric Galton-Watson tree (with parameter $1/2$). The root vertex is the first vertex of the spine. 
\begin{figure}[!h]
\begin{center}
\includegraphics[height=4cm]{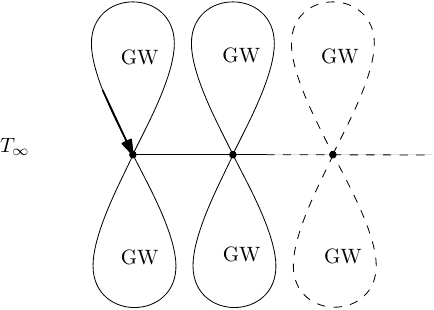}
\caption{ \label{tinftyy} An illustration of $T_{\infty}$.}
\end{center}
\end{figure}
\begin{rek} \label{size-baised}There exists another equivalent way to define $T_\infty$, see \cite{BK06,LPP95b,LP10}: In this description the vertices of the spine of $T_\infty$ have an offspring distribution which is given by a size-biased version of the geometric distribution of parameter $1/2$, whereas all the other vertices have the standard geometric offspring distribution of parameter $1/2$. This construction has the advantage to be easily extended to a general critical offspring distribution. To recover the construction presented above from this one, note that if $G_{1}$ and $G_{2}$ are two independent geometric variables of parameter $1/2$ then $G_{1}+G_{2}+1$ has the same law as a geometric variable of parameter $1/2$ biased by its size.  
\end{rek}
\subsection{The Mass Transport Principle} 
Before presenting the Mass Transport Principle, let us introduce some notation. A graph $G=(\op{V}(G),\op{E}(G))$ is a pair of sets, $\op{V}(G)$ representing the set of vertices and $\op{E}(G)$ the set of (unoriented) edges. In the following, all the graphs considered are countable, connected, locally finite and simple (no loop or multiple edges). For any pair $x,y \in G$, the \emph{graph distance} $\op{d}_{\op{gr}}^{G}(x,y)$ is the minimal length of a path joining $x$ and $y$ in $G$. For every $r \in \mathbb{Z}_{+}$, the ball of radius $r$ around $x$ in $G$ is the subgraph of $G$ spanned by the vertices at distance less than or equal to $r$ from $x$ in $G$, and is denoted by $B_{G}(x,r)$. A \emph{rooted graph} is a pair $(G,\rho)$ where $\rho\in \op{V}(G)$ is called the root vertex. We will identify two rooted graphs if there is a graph isomorphism between them that maps their roots. The set of rooted graphs can be endowed with a metric, see \cite{BC10b,BS01}.

The Mass Transport Principle was introduced by H\"aggstr\"om to study percolation and was further developed by Benjamini, Lyons, Peres and Schramm \cite{BLPS99}. It is extensively studied in \cite{AL07}. We give here an informal definition and refer to \cite{AL07,BC10b} for more details. A random rooted graph $(G,\rho)$ satisfies the Mass Transport Principle if for every positive measurable function $F(g,x,y)$ that associates with each graph $g$ given with  two distinguished vertices $x$ and $y$ of $g$ an ``amount of mass'' sent from $x$ to $y$ in $g$, we have 
\begin{eqnarray} \mathbb{E}\left[\sum_{x \in \mathrm{V}(G)} F(G,\rho,x) \right ] &=&  \mathbb{E}\left[\sum_{x\in \mathrm{V}(G)} F(G,x,\rho)\right] \label{MTP}. \end{eqnarray}
Such a random graph is called \emph{unimodular}.

\begin{corollary}\label{mtp} Let $(\widetilde{T}_{\infty},\tilde{\rho})$ be the random infinite rooted tree obtained from $(T_\infty,\rho)$ after biasing by the inverse of the degree of the root vertex $\rho$, that is 
 \begin{eqnarray*} \mathbb{E}\left[ f(\widetilde{T}_\infty,\tilde{\rho})\right] &=& \frac{\mathbb{E}\left[ \mathrm{deg}(\rho)^{-1}f(T_\infty,\rho)\right]}{\mathbb{E}\left[ \mathrm{deg}(\rho)^{-1}\right]}, \end{eqnarray*} for any positive Borel function $f$. Then  $(\widetilde{T}_{\infty},\tilde{\rho})$  obeys the Mass Transport Principle.
\end{corollary} 
\proof By \cite[Proposition 2.5]{BC10b} it is enough to check that the random rooted graph $(T_{\infty},\rho)$ is reversible in the sense of \cite[Definition 1]{BC10b}. This easily follows from equation \eqref{reversible} and Lemma \ref{convergence}. 
\endproof

\subsection{Volume estimates} 
\label{geometricestimates}Let us introduce some notation.  Recall the construction of the tree $T_\infty$. For $i,j\geq 0$ denote $L_{i}(j)$ (resp.\ $R_{i}(j)$) the number of vertices at the $j$-th generation in the tree grafted on the left (resp.\,right) of the $i$-th vertex of the spine of $T_{\infty}$. It is easy to see that for every $i \geq 0$ the process $(L_{i}(j))_{j\geq 0}$ is a martingale in its own filtration. By standard calculations on geometric distributions we have
\begin{eqnarray}
 \mathbb{E}[ L_{i}(j)] &=& 1 \label{star}  \\
 \mathbb{E}\left[ ( L_{i}(j))^2\right] &=&1 +2j. \nonumber \end{eqnarray}
 Furthermore, by the martingale property of $(L_{i}(j))_{j\geq0}$,  $\mathbb{E}[L_{i}(j)L_{i'}(j')]$ is equal to $1$ if $i\ne i'$ and equals $1+2j$ if $i=i'$ and $j\leq j'$. In particular we can get estimates about the volume of the ball $B_{T_{\infty}}(\rho,r)$ of radius $r$ around $\rho$ in $T_{\infty}$, 
 \begin{eqnarray}
\quad \# B_{T_{\infty}}(\rho,r) &=& r+\sum_{i=0}^{r-1} \sum_{j=1}^{r-i} \big(L_{i}(j) + R_{i}(j)\big), \nonumber \\
\quad  \mathbb{E}[\# B_{T_{\infty}}(\rho,r)] & \sim & r^2.\label{mass} \\
 \mathbb{E}\big[(\# B_{T_{\infty}}(\rho,r))^2\big] & \sim &  \frac{7r^4}{6}.\label{var}
 \end{eqnarray}
 Furthermore, for every $i,h \in \{0,1,2,\ldots\}$ the sum $\sum_{j=0}^{h}L_{i}(j)$ is the size of a critical geometric Galton-Watson tree cut at height $h$. By classical results on critical Galton-Watson trees with finite variance (see e.g. \cite{AN72}), for every $A>0$, there exists $\varepsilon>0$ such that
 $$ \liminf_{h \to \infty} h \cdot \mathbb{P}\left( \sum_{j=0}^{h} L_{i}(j) > \varepsilon h^2\right)  \geq A.$$
 It follows from the last display that the probability that one of the $\lfloor r/2\rfloor$ first trees grafted on the left-hand side of the spine has a size larger than $\varepsilon r^2/4$ is asymptotically at least $1-e^{-A}$ as $r \to \infty$. In particular on this event the number of vertices of $B_{T_{\infty}}(\rho,r)$ is larger than $\varepsilon r^2/4$. Combining this argument with a Markov inequality using equation \eqref{mass} 
 we deduce that 
  \begin{eqnarray}
\lim_{\lambda \to \infty} \inf_{r \geq 0} \mathbb{P}\big( \lambda^{-1} r^2\leq\# B_{T_{\infty}}(\rho,r) \leq \lambda  r^2\big) & =&1. \label{good}
 \end{eqnarray}

 \section{Proof of Theorem \ref{main}}

Let us define precisely what we mean by the simple random walk over $\Z^d$ indexed by a random infinite rooted tree $(T,\rho)$:  Conditionally on $T$ we assign to each edge of $T$ an independent variable uniformly distributed over the symmetric set of the standard basis elements and their inverses of $\mathbb{Z}^d$. For every vertex $u$ in $T$, the sum  of the assigned vectors along the only geodesic path from the root $\rho$ to the vertex $u$ is denoted by $\mathcal{S}_{ ({T},\rho)}(u)$ (note that $ \mathcal{S}_ {(T,\rho)}$ depends on $T$ and on $\rho$). 

 This defines a random function $\mathcal{S}_{ (T,\rho)} : T \to \mathbb{Z}^d$ from the vertices of $T$ to the vertices of $\Z^d$ such that $ \mathcal{S}_ {(T,\rho)}(\rho)=0$. When $ (T,\rho)=(T_{\infty},\rho)$ this function is called \emph{the critical simple random walk snake on $\Z^d$}. We say that the snake is recurrent if almost surely infinitely many vertices of the tree are mapped to $0$ (the origin of $\Z^d$) $i.e.$ 
$$ \# \mathcal{S}_{ ({T}_{\infty},\rho)}^{-1}(\{0\}) = \infty, \quad {\mbox{a.s.}},$$ where the almost surely is in the big probability space of trees and assignments. The snake is transient if $\# \mathcal{S}_{ ({T}_{\infty},\rho)}^{-1}(\{0\}) < \infty$ almost surely.

\subsection{Transience for $d\geq 5$}
For every $r \geq 1$, we denote the set of vertices of $T_\infty$ at distance exactly $r$ from the root vertex $\rho$ by $\partial B_{T_\infty}(\rho,r)$. With the notation introduced in Section \ref{geometricestimates} the number of vertices of $\partial B_{T_\infty}(\rho,r)$ is 
\begin{eqnarray}  \#\partial B_{T_\infty}(\rho,r) = 1 + \sum_{i=0}^{r-1} (L_i(r-i)+R_i(r-i)). \label{estdim5}\end{eqnarray}
Conditionally on $T_\infty$, for every $u \in \partial B_{T_\infty}(\rho,r)$, the probability that $\mathcal{S}_{ ({T}_\infty,\rho)}(u) =0$ is the probability that a simple random walk on $\mathbb{Z}^d$ returns to the origin in exactly $r$ steps, which is less than $\kappa r^{-d/2}$ for some $\kappa >0$. Thus, for every $r\geq 1$ the probability that there exists $u \in \partial B_{T_\infty}(\rho,r)$ such that $ \mathcal{S}_{ ({T}_\infty,\rho)}(u) = 0$ is bounded above by 
\begin{eqnarray*}
\mathbb{P}\left( \exists u \in \partial B_{T_\infty}(\rho,r) : \mathcal{S}_{ ({T}_\infty,\rho)}(u) = 0\right) & \leq & \mathbb{E}\left[ \sum_{u \in \partial B_{T_\infty}(\rho,r)} \mathbf{1}_{\mathcal{S}_{ ({T}_\infty,\rho)}(u) = 0}\right] \\
& \leq & \kappa r^{-d/2}\mathbb{E}\left[\# \partial B_{T_\infty}(\rho,r)\right]\\
& \leq & 3\kappa r^{-d/2+1},
\end{eqnarray*}
where we used \eqref{estdim5} and \eqref{star} to compute $\mathbb{E}\left[\# \partial B_{T_\infty}(\rho,r)\right]$. Consequently, for $d\geq 5$ the preceding bound is summable over $r \in \mathbb{Z}^+$, and an application of Borel-Cantelli's lemma shows that $\mathcal{S}_{ ({T}_\infty,\rho)}$ is transient.
\subsection{Recurrence for $d \in \{1,2,3,4\}$} 
Since the path indexed by the infinite line in $T_\infty$ (the one starting from $\rho$) is distributed as a simple symmetric random walk, the snake is obviously recurrent when $ d \leq 2$. We now suppose $ d \in \{3,4\}$. \emph{We will argue by contradiction} and assume that the simple random walk snake on $\mathbb{Z}^d$ with $d \in \{1,2,3,4\}$ supported by ${T}_{\infty}$ is not recurrent, in particular it is easy to see that we have   \begin{eqnarray}\mathbb{P}\left( \mathcal{S}_{ ({T}_{\infty},\rho)}^{-1}(\{0\})= \{\rho\} \right) &>&0,  \label{assump} \end{eqnarray} in words, the root vertex $\rho$ of $T_\infty$ is the only vertex being mapped to the origin of $ \mathbb{Z}^d$ by the snake $ \mathcal{S}_{ ({T}_\infty,\rho)}$ with positive probability. The following key lemma then states that under this assumption, the range of the snake is somehow linear. This will be the only place where we use the Mass Transport Principle \eqref{MTP}.

\begin{lemma} \label{linear} Assume \eqref{assump}, 
 then there exists $c > 0$ such that for every $r\geq0$,
$$ \mathbb{P}\left( \# \mathcal{S}_{(T_{\infty}, \rho)}\big(B_{T_{\infty}}(\rho,r)\big) > cr^2 \right) > c.$$
\end{lemma}

\proof
Consider the tree $(\widetilde{T}_{\infty},\tilde{\rho})$ obtained after biasing $(T_\infty,\rho)$ by the inverse of the degree of $\rho$. By Corollary \ref{mtp}, this random rooted graph satisfies the Mass Transport Principle \eqref{MTP}. Since we have 
 \begin{eqnarray*} \mathbb{P}\left(\# \mathcal{S}_{(T_{\infty},\rho)}\big(B_{T_{\infty}}(\rho,r)\big) > cr^2\right) &\geq& \mathbb{E}\left[\mathrm{deg}(\rho)^{-1}\mathbf{1}_{\# \mathcal{S}_{(T_{\infty}, \rho)}(B_{T_{\infty}}(\rho,r)) > cr^2}\right]\\ & =& \mathbb{E}[ \mathrm{deg}(\rho)^{-1}]\mathbb{P}\left( \# \mathcal{S}_{(\widetilde{T}_{\infty}, \tilde{\rho})}\big(B_{\widetilde{T}_{\infty}}(\tilde{\rho},r)\big)>cr^2\right), \end{eqnarray*} it is enough to prove the lemma for $\widetilde{T}_{\infty}$ instead of $T_{\infty}$. For every rooted tree $(\tau,\rho)$, we denote  by $\Psi(\tau,\rho)$ the probability that the simple random snake on $\mathbb{Z}^d$ supported on $(\tau,\rho)$  reaches $0$ only at $\rho$, that is
\begin{eqnarray*} \Psi(\tau,\rho) &=& \mathbb{P}(\mathcal{S}_{(\tau, \rho)}^{-1}(\{0\})= \{\rho\}). \end{eqnarray*}
 The function $\Psi$ is thus a positive Borel function over the set of all rooted trees. Notice that if $(\tau,\rho)$ is fixed and if $u \in \tau$ then $ \mathcal{S}_{(\tau,\rho)}- \mathcal{S}_{(\tau,u)}$ has the same distribution as $ \mathcal{S}_{(\tau,u)}$. Thus $\Psi(\tau,u)$ is the probability that the snake $ \mathcal{S}_{(\tau,\rho)}$ is one-to-one at the point $u$ that is   \begin{eqnarray*}  \Psi(\tau,u) & =& \mathbb{P} \Big(\mathcal{S}_{(\tau,\rho)}^{-1}\big( \{ \mathcal{S}_{(\tau,\rho)}(u)\}\big) = \{u\}\Big).  \end{eqnarray*} Using the Mass Transport Principle \eqref{MTP} on $(\widetilde{T}_\infty,\tilde{\rho})$ with the function 
$$ F(G,x,y) = \Psi(G,x) \mathbf{1}_{\mathrm{d}_{\mathrm{gr}}^G(x,y) \leq r},$$
we get 
\begin{eqnarray*} \displaystyle \mathbb{E}\left[\Psi(\widetilde{T}_\infty,\tilde{\rho}) \#B_{\widetilde{T}_\infty}(\tilde{\rho},r) \right] \quad =& \displaystyle \mathbb{E}\left[ \sum_{x \in B_{\widetilde{T}_\infty}(\tilde{\rho},r)} \Psi(\widetilde{T}_\infty,x) \right]\nonumber \\
= & \displaystyle \mathbb{E}\left[ \sum_{x \in B_{\widetilde{T}_\infty}(\tilde{\rho},r)} \mathbf{1}_{\mathcal{S}_{(\widetilde{T}_{\infty}, \tilde{\rho})}^{-1}( \{\mathcal{S}_{(\widetilde{T}_\infty, \tilde{\rho})}(x)\})=\{x\}}\right] & = \ \  \mathbb{E}\Big[ \#I_r\Big],  \end{eqnarray*}   where $I_r$ is the set of vertices in $B_{\widetilde{T}_\infty}(\tilde{\rho},r)$ at which the snake $ \mathcal{S}_{( \widetilde{T}_\infty, \tilde{\rho})}$ is one-to-one. We obviously have $ \# I_r = \#\mathcal{S}_{( \widetilde{T}_\infty, \tilde{\rho})} (I_r)  \leq \#\mathcal{S}_{( \widetilde{T}_\infty, \tilde{\rho})} (B_{ \widetilde{T}_\infty}(\tilde{\rho},r))$ yielding to 
 \begin{eqnarray} \displaystyle \mathbb{E}\left[\Psi(\widetilde{T}_\infty,\tilde{\rho}) \#B_{\widetilde{T}_\infty}(\tilde{\rho},r) \right] & \leq &  \mathbb{E}\Big[\#\mathcal{S}_{( \widetilde{T}_\infty, \tilde{\rho})} (B_{ \widetilde{T}_\infty}(\tilde{\rho},r)) \Big]. \label{MTPT}  \end{eqnarray}
By definition of $\widetilde{T}_{\infty}$, \eqref{good} still holds if we replace $T_{\infty}$ by $\widetilde{T}_{\infty}$. Similarly our condition \eqref{assump} which states that $\Psi({T}_\infty,\rho) >0$ with positive probability implies $\Psi(\widetilde{T}_\infty,\tilde{\rho}) >0$ with positive probability as well. Using these remarks we deduce that the left hand side in \eqref{MTPT} is always larger than some constant times $r^2$: 
 \begin{eqnarray} \inf_{r\geq 1} r^{-2}\mathbb{E}\left[\Psi(\widetilde{T}_\infty,\tilde{\rho}) \#B_{\widetilde{T}_\infty}(\tilde{\rho},r) \right] > 0. \label{lhs}   \end{eqnarray}Let us turn to the right-hand side $\mathbb{E}[ \# \mathcal{S}_{(\widetilde{T}_{\infty},\tilde{\rho})}(B_{\widetilde{T}_{\infty}}(\tilde{\rho},r))]$. Note first that we have $$\# \mathcal{S}_{(\widetilde{T}_{\infty},\tilde{\rho})}(B_{\widetilde{T}_{\infty}}(\tilde{\rho},r)) \leq \#B_{\widetilde{T}_{\infty}}(\tilde{\rho},r).$$ Let $\lambda>0$ and set  $f(x) = x \mathbf{1}_{x > \lambda r^2}$, we get  
\begin{eqnarray}
\mathbb{E}\left[ f\Big(\# \mathcal{S}_{(\widetilde{T}_{\infty},\tilde{\rho})}\big(B_{\widetilde{T}_{\infty}}(\tilde{\rho},r)\big)\Big) \right]  \leq \mathbb{E}\left[ f\big(\#B_{\widetilde{T}_{\infty}}(\tilde{\rho},r)\big) \right] 
  \leq \mathbb{E}[ \mathrm{deg}(\rho)^{-1}]^{-1}\mathbb{E}\left[ f\big(\#B_{T_{\infty}}(\rho,r)\big)  \right]. \nonumber 
  \end{eqnarray} Applying Cauchy-Schwarz inequality we obtain thanks to \eqref{mass} and \eqref{var}
   \begin{eqnarray*} \mathbb{E}\left[ f\big(\#B_{\widetilde{T}_{\infty}}(\tilde{\rho},r)\big)  \right] & \leq & \Big(\mathbb{E}\left[\left(\#B_{T_{\infty}}(\rho,r)\right)^2\right] \mathbb{P}\left( \#B_{T_{\infty}}(\rho,r) > \lambda r^2 \right) \Big) ^{1/2} \leq C r^2\lambda^{-1/2},  \end{eqnarray*}
for some positive constant $C>0$.   Hence 
    \begin{eqnarray} \sup_{r >0} r^{-2}\mathbb{E}\left[ \# \mathcal{S}_{(\widetilde{T}_{\infty},\tilde{\rho})}\big(B_{\widetilde{T}_{\infty}}(\tilde{\rho},r)\big) \mathbf{1}_{\# \mathcal{S}_{(\widetilde{T}_{\infty},\tilde{\rho})}(B_{\widetilde{T}_{\infty}}(\tilde{\rho},r))> \lambda r^2} \right] &\xrightarrow[\lambda \to\infty]{}& 0.  \label{pasgrand} \end{eqnarray}
Combining \eqref{MTPT}, \eqref{lhs} and \eqref{pasgrand}, we deduce that there exists a constant $c>0$ such that 
$$ \mathbb{P}(\# \mathcal{S}_{(\widetilde{T}_{\infty},\tilde{\rho})}(B_{\widetilde{T}_{\infty}}(\tilde{\rho},r)) > cr^2) > c,$$ which is the desired result.\endproof

Now consider $\eta>0$ small enough so that 
$$ A_{\eta,r}(T_{\infty}) = \left\{ u \in B_{T_{\infty}}(\rho,r) \mbox{ in a tree grafted on the spine of }T_{\infty} \mbox{ before } \lfloor \eta r \rfloor\right\} $$ 
has a cardinal less than $cr^2/2$ with probability at least  than $1-c/2$, independently of $r \geq 1$. Then by Lemma \ref{linear}, for every $r \geq1$,  the set $B_{T_{\infty}}(\rho,r) \backslash A_{\eta,r}(T_{\infty})$ which consists of the vertices within distance $r$ of $\rho$ that are linked to the spine after the $\lfloor\eta r\rfloor$-th vertex has an image by $\mathcal{S}_{(T_{\infty},\rho)}$ of size larger than $c r^2/2$ with probability at least $c/2$,
$$ \mathbb{P}\left( \#\mathcal{S}_{(T_{\infty},\rho)}\big( B_{T_{\infty}}(\rho,r) \backslash A_{\eta,r}(T_{\infty})\big) > cr^2/2\right) > c/2.$$
By diffusivity bounds on the simple random snake on $T_{\infty}$ (see \cite{CMM10} for the case $d=1$ which is easily extended to $d\geq 2$, see also \cite{Kes95}), one can find $M>0$ large enough so that the image of $B_{T_{\infty}}(\rho,r)$ by the snake is contained in $B_{\Z^d}(0,M \sqrt{r})$ with probability larger than $1-c/4$. If $d \in \{1,2,3\}$ we already reached a contradiction since $\#B_{\Z^d}(0,M \sqrt{r})  \leq 8M^3 r^{3/2}$ and thus $B_{\Z^d}(0,M \sqrt{r}) $ cannot contain a set of size of order $r^2$ with positive probability for $r$ large enough. 

We now suppose $d=4$. Summing-up, for every $r \geq1$, with a probability at least $c/4$, the image of $B_{T_{\infty}}(\rho,r) \backslash A_{\eta,r}(T_{\infty})$ by $\mathcal{S}_{(T_{\infty},\rho)}$ is a random set in $ \mathbb{Z}^4$ composed of at least $cr^2/2$ different vertices and whose diameter is less than $M \sqrt{r}$. 

\begin{lemma} The point $0$ (origin of $ \mathbb{Z}^4$) is in the image of $B_{T_{\infty}}(\rho,r) \backslash A_{\eta,r}(T_{\infty})$ by $\mathcal{S}_{(T_{\infty},\rho)}$ with a probability bounded away from $0$ independently of $r \geq1$.
\end{lemma}

\proof All the vertices of $B_{T_{\infty}}(\rho,r) \backslash A_{\eta,r}(T_{\infty})$ are linked to $\rho$ by the same first $\lfloor \eta r \rfloor$ vertices of the spine of $T_{\infty}$. Besides, the increments along the first $\lfloor \eta r \rfloor$ edges of the spine are independent  of the structure of $B_{T_{\infty}}(\rho,r) \backslash A_{\eta,r}(T_{\infty})$ and also independent of the increments of the snake along the edges of $B_{T_{\infty}}(\rho,r) \backslash A_{\eta,r}(T_{\infty})$. 

Denote by $S$ the image of  $B_{T_{\infty}}(\rho,r) \backslash A_{\eta,r}(T_{\infty})$ by $\mathcal{S}_{(T_{\infty},\rho)}$ translated such that the image of the $\lfloor \eta r \rfloor$-th  vertex of the spine is  $0$. Let $E_{r}$ be the event that $S$ is of size larger than $c r^2/2$ and of diameter less than $M \sqrt{r}$.  By the arguments developed before the lemma, we have $ \mathbb{P}(E_r) > c/4$ for every $r \geq 1$. Besides, the image of $B_{T_{\infty}}(\rho,r) \backslash A_{\eta,r}(T_\infty)$ by $\mathcal{S}_{(T_{\infty},\rho)}$ is a random translation of the set $S$ by an independent $\lfloor\eta r\rfloor$-steps random $X_{\lfloor\eta r\rfloor}$ walk on $\Z^4$. 
We claim that the conditional probability $\mathbb{P}( \mathbf{1}_{0 \in S+ X_{\lfloor \eta r\rfloor }} | E_{r})$ is bounded away from $0$ independently of $r \geq 1$. Indeed conditionally on $S$ and $E_{r}$, the variable $$I=\sum_{x \in S+X_{\lfloor \eta r\rfloor }} \mathbf{1}_{x =0}$$ takes values in $\{0,1\}$ and its expectation is 
$$ \mathbb{E}[I \mid E_r] =  \mathbb{E}\left[\left. \sum_{x \in S} \mathbb{P}(X_{\lfloor \eta r\rfloor } = -x)\,  \right| E_r \right] \geq \frac{cr^2}{2} \cdot \kappa r^{-2},$$   for some constant $\kappa>0$ depending on $\eta$ and $M$ only. Hence, on the event $E_r$, the expectation of $I$ is bounded away from $0$, thus $\mathbb{P}(I>0\mid E_{r}) = \mathbb{E}[I\mid E_{r}]$ is bounded away from $0$. Since $\mathbb{P}(E_{r}) \geq c/4$ for all $r\geq 0$, the lemma is proved.
\proof For $\eta >0$ fixed, denote $ \mathcal{B}_r = B_{T_{\infty}}(\rho,r) \backslash A_{\eta,r}(T_{\infty})$ to simplify notation. Thanks to the preceding lemma, the probability that $0$ belongs to the image of $ \mathcal{B}_r$ by the snake is bounded from below by some positive constant $c>0$ independent of $r \geq 0$. On the other hand we have,   \begin{eqnarray*} \mathbb{P} \Big( 0 \in \mathcal{S}_{(T_\infty,\rho)}(  \mathcal{B}_r)\cap \mathcal{S}_{(T_\infty,\rho)}(  \mathcal{B}_{r'}) \Big) & \xrightarrow[r'\to\infty]{} & \mathbb{P}\Big( 0 \in \mathcal{S}_{(T_\infty,\rho)}(  \mathcal{B}_r)\Big) \mathbb{P} \Big( 0 \in \mathcal{S}_{(T_\infty,\rho)}(  \mathcal{B}_{r'})\Big).  \end{eqnarray*}
In words, the events $\{ 0 \in \mathcal{S}_{(T_\infty,\rho)}(  \mathcal{B}_r)\}$ and $\{0 \in \mathcal{S}_{(T_\infty,\rho)}(  \mathcal{B}_{r'})\}$ are asymptotically independent as $r' \to \infty$.  Hence, by an adaptation of Borel-Cantelli's lemma we deduce that $\{ 0 \in \mathcal{S}_{(T_\infty,\rho)}(  \mathcal{B}_r)\}$ occurs for infinitely many $r$'s. This implies recurrence, contradiction.\endproof

\section{Extensions and comments} \label{extensions}
\subsection{General trees}
Theorem \ref{main} still holds for more general critical Galton-Watson trees conditioned to survive. Namely, if $\xi$ is a critical offspring distribution with finite variance, Lemma \ref{convergence} is still true and the limiting tree called the $\xi$-Galton-Watson tree conditioned to survive can be described as in Remark \ref{size-baised}. The geometric estimates of Section \ref{geometricestimates} can be adapted for this infinite tree. However the analogue of Proposition \ref{prop:rev} and its Corollary \ref{mtp} are not true in this general setting. To bypass this difficulty, one can rely on the trick introduced in \cite{LPP95} and use the so-called \emph{augmented Galton-Watson} measure. To be precise, let  $T^\xi$ and $T^{\xi}_{\infty}$ be respectively a critical $\xi$-Galton-Watson tree and a critical $\xi$-Galton-Watson tree conditioned to survive. We suppose that $T^{\xi}$ and $T_{\infty}^{\xi}$ are independent and we define the tree $\overline{T}_{\infty}^\xi$ obtained by joining the roots of the trees $T^{\xi}$ and $T^{\xi}_{\infty}$ by an edge. The root $\rho$ of this tree is with probability $1/2$ the root vertex of $T^\xi$ and with probability $1/2$ the root vertex of $T_{\infty}^\xi$.

\begin{figure}[!h]
 \begin{center}
 \includegraphics[height=4cm]{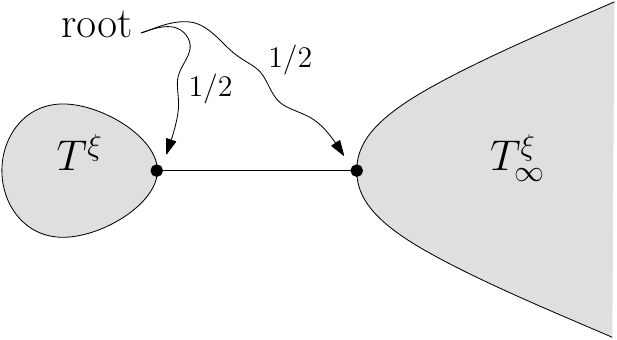}
 \caption{An illustration of the augmented critical Galton-Watson tree conditioned to survive.}
 \end{center}
 \end{figure}
 
 Then one can show that $(\overline{T}_{\infty}^\xi,\rho)$ is stationary and reversible (in the sense of \cite[Definition 1]{BC10b}) thus once biased by the inverse of the degree of the root it satisfies the Mass Transport Principle, see \cite{LPP95,LP10} for closely related proofs.

\subsection{Alternative proof and extension}
It is likely that a proof of Theorem \ref{main} would also follow from an adaptation of estimates for superBrownian motion to discrete snakes (see \cite{DIP89,LG99}), for example by showing
$$ \mathbb{P}\big( \mathcal{S}_{(T,\rho)} \mbox{ reaches } x \in \mathbb{Z}^4\big) \asymp \frac{1}{\log(|x|)|x|^2},$$ where the snake $ \mathcal{S}_{(T,\rho)}$ runs over a critical geometric Galton-Watson tree $(T,\rho)$ \emph{not} conditioned to be infinite. This approach has been carried out in an unpublished work of H. Kesten and Y. Peres (personal communication). In particular such estimates would be needed to answer the following questions:
 \begin{quest}
 What is the variance of the number of returns to $0$ by the snake $ \mathcal{S}_{(T_\infty,\rho)}$ restricted to the ball of radius $r$ in $T_{\infty}$? What is the range of the snake restricted to the ball of radius $r$ in $T_{\infty}$?
\end{quest}

  However we believe that our more abstract argument relying on the  Mass Transport Principle and rough volume estimates could be used in a more general setting including e.g.\,proving that a simple random snake on $\mathbb{Z}^{3}$ indexed by a critical Galton-Watson conditioned to be infinite whose offspring distribution is in the domain of attraction of a stable law of parameter $3/2$ is recurrent. We end this note with the following (possibly hard) question, generalizing intersecting probabilities for SRW in two dimensions:
  \begin{quest}
  Consider two independent snakes $\mathcal{S}_{(T_{\infty},\rho)}$ and $\mathcal{S}'_{(T_{\infty}',\rho')}$ on $ \mathbb{Z}^4$ and starting from $0$ (origin of $ \mathbb{Z}^4$) at the roots $\rho,\rho'$ of $T_{\infty},T'_{\infty}$. For $r \geq 0$, estimate the probability that the two snakes supported by the balls of radius $r$ in the two trees intersect only at $0$, that is estimate
  $$ \mathbb{P}\left( \mathcal{S}_{(T_{\infty},\rho)}^{-1}(B_{T_{\infty}}(\rho,r)) \cap {\mathcal{S}^{'-1}_{(T'_{\infty},\rho')}}(B_{T'_{\infty}}(\rho',r)) = \{0\} \right).$$
\end{quest}

\clearpage

\bibliographystyle{abbrv}

\end{document}